\documentclass{amsart}
\usepackage{amscd,amssymb}

\theoremstyle{plain}
\newtheorem{thm}{Theorem}[section]
\newtheorem{lem}[thm]{Lemma}
\newtheorem{definicja}[thm]{Definition}
\newtheorem{cor}[thm]{Corollary}
\newtheorem{prop}[thm]{Proposition}

\theoremstyle{definition}

\theoremstyle{remark}

\numberwithin{equation}{section}

\begin{document}
\title{Isometry types of profinite groups}
\author{A.Ivanov (Wroc{\l}aw)}
\maketitle

\bigskip

{\em Abstract.} 
\footnote{{\bf Key words and phraces:} rooted trees, isometry groups, 
profinite groups, conjugacy classes, interpretability} 
\footnote{the research was supported by KBN grant 1 P03A 025 28}
Let $T$ be a rooted tree and $Iso(T)$ be 
the group of its isometries. 
We study closed subgroups $G$ of $Iso(T)$ with respect 
to the number of conjugacy classes of $Iso(T)$ having 
representatives in $G$. 

\bigskip

\section{Introduction}

Let $T$ be an infinite locally finite tree with the root
$\emptyset$. 
By $layer_n$ we denote the set of all elements of $T$, 
$n$-distant from $\emptyset$ (i.e. elements of level $n$). 
We always assume that there is an infinite sequence 
$\bar{m}=m_1 ,m_2 ,...$ of natural numbers such that 
$T$ is of the form $T_{\bar{m}}$, the tree where 
all elements of $layer_n$ are of valency $m_{n+1}+1$ for $n>0$. 
Let $T_n = \bigcup _{i\le n} layer_i$ and 
$T_{\langle n \rangle}$ be the tree of the isomorphism type 
of $T_w$, the tree hanging from $w\in layer_n$.  \parskip0pt 

The group $Iso(T)$ of all isometries of $T$ fixing $\emptyset$ 
is a profinite group with respect to canonical homomorphisms
$\pi_{n}: Iso(T)\rightarrow Iso(T_{n})$. 
In fact the following examples represent  
all the main fields in the subject. 

Consider the group $GL_d (\mathbb{Z}_p )$ and 
its action on $D=(\mathbb{Z}_p )^d$. 
Let $T_{ad}$ be the tree of all cosets in $(\mathbb{Z}_p )^d$ 
with respect to all subgroups of the form $p^k (\mathbb{Z}_p )^d$.  
If $G$ is a closed subgroup of $GL_{d}(\mathbb{Z}_p )$,  
then $D$ becomes a continuous $G \mathbb{Z}_p$-module. 
Thus $G$ acts by isometries on $T_{ad}$ and fixes the root $D$. 
A very close example is the action of $PGL_2 (\mathbb{Z}_p )$ 
on the tree of lattices in $\mathbb{Q}_p \times \mathbb{Q}_p$ 
(see \cite{Serre}). 
\parskip0pt 

Completely different groups are provided by Grigorchuk-type 
constructions. 
In \cite{Grig} {\em branch groups} are introduced, 
which in some sense can be considered as an axiomatization 
of well-known examples (their profinite completions) 
of Burnside groups found by Grigorchuk.  
  
Are there any global (=model-theoretic) properties 
giving interesting dividing lines in these examples ?  

In this paper we concentrate on characterizations  
concerning diversity of elements occurring in the group. 
Our basic property can be described as follows. 
Consider the relation $\sim$ of conjugacy in the group
$Iso(T)$ of all isometries of $T$. 
It is easy to see that every conjugacy class is a closed 
subset of $Iso(T)$.
Let $G$ be a closed subgroup of $Iso(T)$. 
We say that $G$ has {\em a small (resp. countable) number 
of isometry types} if $G$ meets only $<2^{\aleph_{0}}$ 
(countably many) conjugacy classes of $Iso(T)$. 
Since $G/\sim$ is analytic, smallness
is equivalent to countability. \parskip0pt 

We study profinite groups which can be embedded 
into $Iso(T)$ as closed subgroups with a small 
number of isometry types. 
It seems to us that this class natuarally arises 
in many situations. 
One of them has model-theoretic flavor and is worth 
mentioning here. 
For closed $G<Iso(T)$ consider the corresponding 
inverse system 
$G_0 \leftarrow G_1 \leftarrow ... \leftarrow G_i \leftarrow ...$ 
with $G_i\le Iso(T_i )$. 
Any automorphism of $G$ fixing all $G_i$ (as sets) 
is called a {\em profinite automorphism} of $G$. 
By $Aut^* (G)$ we denote the group of all profinite 
automomorphisms. \parskip0pt 
   
Assume that any profinite automorphism of $G$ is
induced by an isometry of $T$ as an inner automorphism of
$Iso(T)$ (by Proposition 8.1 from \cite{LN} this happens if
for any $w \in T$ the point-wise stabilizer of
$T\setminus T_{w}$ is level-transitive on $T_{w}$).
If the {\em number of $Aut^{*}(G)$-orbits on $G$ is at most countable}
then $G$ has a small number of isometry types. 
The former condition has been introduced by Newelski 
in \cite{N} as a counterpart of small theories 
from model theory.
The main (still open) conjecture in these respects 
states that profinite groups with a small number of 
$Aut^{*}(G)$-orbits have open abelian subgroups. 
Some partial confirmation has been obtained by F.Wagner 
in \cite{Wagner}. 
\parskip0pt 

In fact Newelski's conjecture has become the first motivation 
for our interests to isometry groups with a small number 
of isometry types. 
Notice that when $G$ is an inverse limit of 
a system of finite groups $G_0 \leftarrow G_1 \leftarrow ...$ 
$\leftarrow G_i \leftarrow ...$ with kernels $N_i$ of 
the corresponding canonical homomorphisms $G\rightarrow G_i$, 
then $G$ acts by left multiplication on the tree $T_C (G)$ of 
left cosets of all $N_i$, $i\in\omega$ (any edge of the tree 
is defined by a pair of the form $gN_{i+1} \subset gN_{i}$).  
In this case any profinite automorphism of $G$ 
(which stabilises all $N_i$) becomes an isometry of $T_C (G)$. 
Thus if $G$ has a small number of $Aut^{*}(G)$-orbits, then 
$G$ has a small number of isometry types with respect 
to this tree action. 

It has become very surprising for us that the smallness 
condition above has very strong abstract consequences 
concerning the structure of $G$. 
In patricular it turns out that we can say quite a lot 
about $G$ when $G$ is nilpotent.  

We also study the related property of definability of 
the action in $\mathbb{Q}_p$ and show that some natural 
zeta functions arising in these respects are rational. 
We compare our approach with the one suggested by 
Z.Charzidakis in \cite{zoe}. 
We will see that $\omega$-categoricity appearing 
in \cite{zoe} is very far from smallness of the number 
of isometry types (at least in the case of examples above).

\subsection{Preliminaries} 

Let $T=T_{\bar{m}}$, where $\bar{m}$ is an infinite sequence 
of positive natural numbers. 
For every $n$ we define an alphabet $A_n =\{ a^n_1 ,...,a^n_{m_n}\}$ 
and consider $T$ as the union of all finite products 
$A_1 \times ...\times A_n$, $n\in\omega$, interpreted 
as paths of $T$ started at $\emptyset$: 
here we identify $layer_n$ with $A_1 \times ...\times A_{n}$ 

If $v$ is a vertex of $T$ and $g\in Iso(T)$, then by $g(v)$ 
we usually denote the result of $g$ applied to $v$. 
On the other hand we sometimes use the form $v^g$, because 
it becomes more convenient when one considers iterated wreath products.  
Then as in \cite{Grig} any $g\in Iso(T)$ is regarded as 
a labelling $\{ \gamma (v)\}_{v\in T}$ of the vertices by elements 
of the symmetric group which acts on the edges below the vertex. 
The effect of $g\in Iso(T)$ on a vertex $u$ 
corresponding to $p_u =(a_1 ,...,a_n )$, a path from 
$\emptyset$  with $a_i \in A_i$, is given by the formula: 
$$
p^g_u = ( a^{\gamma (\emptyset )}_1 ,a^{\gamma (a_1 )}_2 ,
..., a^{\gamma (a_1 ...a_{n-1})}_n ), 
$$ 
where by $a^{\gamma (v)}$ we denote the result of the permutation 
$\gamma (v)$ applied to the edge $a$.  
The rule for the composing automorphisms is given by 
the following formulas (see \cite{Grig}): 
\begin{quote} 
Let $h=\{\delta (v)\}_{v\in T}$, 
$f=\{\sigma (v)\}_{v\in T}$ and $g=\{ \gamma (v) \}_{v\in T}$. \\ 
If $h=fg$, then $\delta (v) = \sigma (v)\gamma (v^f )$. \\ 
If $f=g^{-1}$, then $\sigma (v) = (\gamma (v^{g^{-1}}))^{-1}$. 
\end{quote} 
  
Now assume that for every $n$ we have chosen a permutation group 
$P_i \le Sym (A_i )$. 
Then the subgroup of all labellings $\{ \gamma (v)\}_{v\in T}$ 
with $\gamma (v)\in P_{i+1}$ for $v \in layer_i$ is called 
{\em the iterated wreath product} of $(P_i ,A_i )$, $i\in \{ 1,2,...\}$ 
and is denoted by $\wr_{\omega} (P_i ,A_i )$.  

With any $v\in T$ (an end $\mathbf{e}\in \partial T$ resp.) and
a group $G\le Iso(T)$ we associate the stabilizer $St_G (v)$ 
($St_G ({\bf e})$ resp.). 
The level stabilizer is defined by 
$St_G (n) = \bigcap \{ St_G (v): |p_v |=n\}$.  
It is clear that $St_G (n)= Ker \pi_n$, where 
$\pi_n :Iso (T) \rightarrow Iso(T_n )$ is defined as above. 
On the other hand $St_{G} (n)$ is naturally identified 
with a subgroup of $\prod \{ Iso (T_v ): v\in layer_n\}$. 

The {\em rigit stabilizer} of $v$ is defined to be 
$Rs_G (v) = \{ g\in Iso(T): g$ fixes $T\setminus T_v$ pointwise$\}$. 
The rigit stabilizer of the $n$-th level $Rs_G (n)$ is defined 
to be the group generated by $\bigcup \{ Rs_G (v): v\in layer_n \}$. 
It is clear that $Rs_G (n)= \oplus \{ Rs_G (v): v\in layer_n \}$. 

We now define (weakly) branch groups. 
They have been introduced by R.Grigorchuk in \cite{Grig}. 
They are one of the main objects in the area. 

\begin{def} \label{branch} 
Let $T$ be a locally finite tree, $T=T_{\bar{m}}$. 
A closed subgroup $G<Iso(T)$ is called a branch group with 
respect to $T$ if \\ 
(i) $G$ acts transitively on each layer of the tree, and \\ 
(ii) for each $n\ge 1$ there exists closed 
$L_n \le Iso(T_{\langle n\rangle })$ 
such that the direct product $H_n = L^{layer_n}_n$ is normal 
and of finite index in $G$ (each factor of the product 
acts on the subtree hanging from the corresponding vertex of 
the $n$-th layer). 
\end{def}  

Condition (ii) of this definition is equivalent to 
$|G:Rs_G (n)|<\infty$ for every $n\in\omega$ \cite{Grig}. 
If we replace (ii) by the condition that $Rs(n)$ is 
always infinite, then we obtain a definition of 
{\em weakly branch groups}.  

Let $g \in Iso(T)$. 
By $T_g$ we denote the {\em orbit tree} of $g$, 
i.e. the set of all $g$-cycles $T/\langle g\rangle$ 
with respect top the natural ajacency. 
The orbit tree is labelled by natural numbers as follows: 
to each orbit $\langle g\rangle v$ we assosiate 
the size $|\langle g\rangle v|$. 

The following theorem from \cite{GNS} will be one of the main 
tools of the paper: 
\begin{quote} 
isometries $g$ and $h$ are conjugated in $Iso(T)$ if and only if 
the corresponding orbit trees $T_g$ and $T_h$ are isomorphic 
by a labelling preserving isomorphism. 
\end{quote} 

It is worth noting that the class of groups having a
faithful representation with a small number of isometry types
is closed under direct sums: if $G_{1}$ and $G_{2}$ have
appropriate realisations on trees of valencies $N_{1}$ and
$N_{2}$, then $G_{1}\oplus G_{2}$ can be realised on a tree
of valency $N_{1}+N_{2}$ so that the number of isometry types
is small.

\subsection{Main results} 

In Section 2 we study profinite groups which can be 
realized as closed subgroups of $Iso(T)$ 
which have a small number of isometry types. 
We show that $\mathbb{Z}_p$ and the direct power 
$H^{\omega}$ of a finite group $H$ have such 
representations in $Iso(T)$ for appropriate $T$. 
We conjecture that nilpotent profinite groups $G$ with 
this property are very close to these examples. 
In Section 2.1 we slightly confirm this showing that 
only finitely many primes can divide $|G|$ and the periodic 
part of $G$ is of finite exponent. 

We show that the smallness condition is so restrictive 
itself that for example weakly branch groups do not have 
a small number of isometry types.  
We will also see here that the tree of 
cosets of kernels of canonical projections naturally 
arises in these respects, and the smallness of the number 
of isometry types implies some interesting properties 
of this tree (Proposition 2.2). 
 
In Section 2.2 we study pronilpotent groups $G$ 
which have a small number of isometry types for 
all faithful continuous representations in isometry 
groups of locally finite trees. 
In the case when the group $G$ is nilpotent we 
show that such a group must be of finite rank. 
This condition becomes sufficient when the 
group is abelian (see Theorem \ref{Z_p} ). 

In Section 3 we study another condition of smallness. 
Let $T$ be a locally finite rooted tree. 
For every natural $n$ we introduce the following 
equivalence relation on $Iso(T)$: 
for $g,g'\in Iso(T)$ we define $g \equiv_{n} g'$ 
if the restrictions of $g$ and $g'$ to $T_n$ are conjugated 
in $Iso(T_n )$.  
\parskip0pt 

For a subgroup $G\le Iso(T)$ we denote by $c^{G}_n$ 
the number of all $\equiv_{n}$-classes meeting $G$. 
In Section 3 we will study the question when 
the zeta function $\Sigma_{n\ge 0} c^{G}_{n} t^n$ 
is rational. 

In general one cannot expect, that this 
happens for a typical $G<Iso(T)$ 
For example Grigorchuk shows in \cite{Grig2} that 
the Hilbert-Poincar\'{e} series of some of his groups 
are not rational.  
But if we smooth up the situation 
(under some model theoretic assumptions) 
we can have very reasonable variants of the question. 
Moreover examples of the introduction with 
$GL_d (\mathbb{Z}_p )$ involved become under our methods. 

The idea looks as follows. 
Let $\Gamma$ be a closed subgroup of $Iso(T)$. 
We say that $g,g'\in Iso(T)$ have the same 
$\Gamma$-type if the exists $h\in \Gamma$ 
such that $g^{h}=g'$.   
It makes sense to study $\Gamma$-types of $Iso(T)$ 
assuming that $\Gamma$ is interpetable in some natural 
model-theoretic objects (not only in the case $\Gamma =Iso(T)$).  
In this situation we say that $G\le Iso(T)$ has a small 
number of $\Gamma$-types if $G$ meets $\le \aleph_0$ 
$\Gamma$-orbits of the conjugacy action of $\Gamma$ 
on $Iso(T)$.  \parskip0pt 
   
For every natural $n$ we introduce the following 
equivalence relation on $Iso(T)$: 
for $g,g'\in Iso(T)$ we define $g \equiv^{\Gamma}_{n} g'$ 
if the restrictions of $g$ and $g'$ to $T_n$ are conjugated by 
an element of the projection of $\Gamma$ to $T_n$.  
Now for the subgroup $G\le Iso(T)$ we denote by $c^{G}_n$ 
the number of all $\equiv^{\Gamma}_{n}$-classes meeting $G$. 
We will really study the question when the zeta function 
$\Sigma_{n\ge 0} c^{G}_{n} t^n$ is rational for some 
$\Gamma$ interpretable in the field $\mathbb{Q}_p$. 
In Section 3 we discuss a number of examples. 
In particular we study a point-stabilizer of the action 
of $SL_2 (\mathbb{Q}_p )$ on the tree of lattices \cite{Serre}. 
In fact our material here is based on a nice paper of E.Hrushovski 
and B.Martin \cite{hrma}. 
\parskip0pt 

We must mention that in the eighties Z.Chatzidakis 
introduced some language for profinite groups and showed 
that model-theoretic investigations of the corresponding 
structures (denoted by $S(G)$) give very strong consequences 
in purely algebraic questions about profinite groups 
and field extensions. 
Moreover structures $S(G)$ arising in this approach, have 
very nice model-theoretic properties. 
In some typical situations they are $\omega$-categorical, 
i.e. small. 
It should be natural to start model-theoretic investigations 
of closed subgroups of $Iso(T)$ with these aspects. 
We fill this gap in Section 4. 
The answer which we obtain here, shows that the majority 
of our examples represent the simplest case of this classification. 
In fact  we show that the dividing line of $\omega$-categoricity of 
structures $S(G)$ is situated inside the class of branch groups. 
We give here some interesting examples. 
In particular the profinite completion of the famous example 
of Grigorchuk is $\omega$-categorical in this sense.    
These results show why the condition of smallness which 
comes from the Newelski's approach is more appropriate 
for our examples.

\section{Groups with a small number of isometry types}

\subsection{Groups having representations with a small number of isometry types}

We now consider groups which can be realized as closed 
subgroups of $Iso(T)$ with a small number of isometry types. 
Let $T$ be a locally finite tree and $G$ be a closed 
subgroup of $Iso(T)$.  
By $\pi_n$ we denote the canonical homomorphism from $G$ to 
the group $G_n$ of all restrictions of $G$ to $T_n$, $n\in \omega$. 

\begin{lem}\label{trans} 
(1) If there is some $n\in\omega$ such that $ker(\pi_n )\setminus \{ 1\}$ 
consists of isometries of the same type, then the group $G$ 
has an open subgroup which is a locally finite pro-$p$-group 
for some prime $p$. 

(2) If there is $n\in\omega$ such that the set of non-diagonal 
pairs of $ker(\pi_n )$ is contained in the same orbit of $Iso(T)$, 
then $ker(\pi_n)$ is isomorphic to the power $(\mathbb{Z}(p))^{\omega}$.  
\end{lem}

{\em Proof.} 
Under the assumptions of the lemma we easily see that 
for every $l\ge n$ all non-trivial elements of the finite group 
$ker(\pi_n )/ker(\pi_l )$ are of the same order. 
Thus all of them are of order $p$, where $p$ is prime. 
This shows that $ker(\pi_n )$ is a pro-$p$-group 
of exponent $p$. 
By a theorem of Zelmanov from \cite{zelmanov} it is 
locally finite. \parskip0pt 

The second statement follows from the fact that every 
$ker(\pi_n )/ker(\pi_l )$ is a nilpotent group and   
has non-trivial centre. 
The assumption says that it must be isomorphic to a finite power 
of $\mathbb{Z}(p)$. 
The rest is obvious. $\Box$

\begin{prop} \label{proptrans} 
Let $T$ be a locally finite tree and $G$ be a closed 
subgroup of $Iso(T)$ with a small number of isometry types.  
Then the following statements hold. 

(1) For every $m\in\omega$ and every $h\in G$ there are some 
$n\in\omega$ and $g\in G\setminus ker(\pi_n )$ 
such that $g\cdot ker(\pi_n)\subseteq h\cdot ker(\pi_m)$ and 
$g\cdot ker(\pi_n)$ consists of isometries of the same type.  

(2) If $g\in G\setminus ker(\pi_n )$, $n\in\omega$, and the set 
of all non-diagonal pairs of $g\cdot ker(\pi_n )$ 
is contained in the same orbit of $Iso(T)$, then $G$ has 
an open subgroup which is a locally finite pro-$p$-group 
for some prime $p$. 
If moreover for every $g_1, g_2 ,g'_1 ,g'_2 \in g\cdot ker(\pi_n )$ 
with $g_1 \not= g_2$, $g'_1 \not= g'_2$ one of the pair 
$(g_1 ,g_2 )$ and $(g'_1 ,g'_2 )$ can be mapped to another 
by an isometry fixing some element of $g\cdot ker(\pi_n )$ 
then $ker(\pi_n )$ is isomorphic to the power 
$(\mathbb{Z}(p))^{\omega}$, where $p$ is prime.  
\end{prop}

{\em Proof.} 
(1) By the Bair category theorem, $Iso(T)$-conjugacy 
classes which are open in $G$ form a dense subset of $G$. 
By the definition of the profinite topology we see the existence 
of $g\cdot ker(\pi_n )$ as in the formulation. 
\parskip0pt 
 
(2) If the set of all non-diagonal pairs of $g\cdot ker(\pi_n )$ 
is contained in the same orbit of $Iso(T)$, then 
each element of $ker(\pi_n )$ can be mapped to any 
other element of $ker(\pi_n )$ by an isometry. 
By Lemma \ref{trans}, $G$ has an open subgroup which 
is a locally finite pro-$p$-group for some prime $p$. 
\parskip0pt 

Assuming the last condition of the formulation consider 
any $ker(\pi_n )/ker(\pi_l )$ with $l\ge n$. 
As we already know this is a group of exponent $p$. 
Let $h_1 \cdot ker(\pi_l)$ belong to the centre of 
$ker(\pi_n )/ker(\pi_l )$, and $h_2 \not= h_3 \in ker(\pi_n )$. 
Let an isometry $\mu$ map $gh_1$ to $gh_2$ and fix the pair $g, gh_3$
(we may assume so). 
Then we have modulo $ker(\pi_l )$: 
$gh_1 \cdot gh_3 = ggh_3 (h_1 )^{g} = g(h_3 )^{g^{-1}}\cdot g(h_1 )^{g}$ 
and 
$gg(h_2 )^{g}h_3 = gh_2 gh_3 =\mu (gh_1 )\mu (gh_3)=\mu (gh_1\cdot gh_3)$ 
$=g(h_3 )^{g^{-1}}g(h_2 )^{g}= ggh_3(h_2 )^{g}$. 
Thus $h_2 ker(\pi_l )$ belongs to the centre of $ker(\pi_n )/ker(\pi_l )$. 
The rest is obvious. $\Box$ 

\bigskip 

\begin{cor} \label{8} 
Let $T$ be a locally finite tree and $G$ be a closed 
subgroup of $Iso(T)$ with a small number of isometry types. 
Then there is no infinite sequence of vertices $v_{1},v_{2},...$
such that $T_{v_{i}}\cap T_{v_{j}}=\emptyset$ for $i\not=j$,
and for every $i$ the point-wise stabilizer of
$T \setminus  T_{v_{i}}$ in $G$ has a non-trivial
action on $T_{v_{i}}$. 
In particular a weakly branch group does not have a small 
number of isometry types.
\end{cor} 

{\em Proof.} 
Assume the contrary. 
Since $T$ is locally finite, we easily see that $G$ does not satisfy 
Proposition \ref{proptrans} (1). 
The rest is obvious. 
$\Box$

\bigskip

{\bf Example.} 
One can build a tree action of the wreath product $\mathbb{Z}(p)wr\mathbb{Z}_p$
with a large number of isometry types.

The following example concerns the case of a tree of valency 
three. 

\bigskip 

{\bf Example.}
Let $T=\{ 0,1,\}^{<{\bf N}}$ be a tree with the root
$\emptyset$ (${\bf N}= \omega\setminus \{ 0\}$; the
elements of the tree are presented by finite sequences).
Let $g\in Iso(T)$ be defined as follows: $g(0) = 1$,
$g(1) = 0$. 
We define $g$ so that at every level $g$
has a unique $q$-cycle: If $v$ consists of $n$ $0$-s
(the beginning of the cycle) then let
$g^{i}(v0)=g^{i}(v)0$, $i < 2^{n}$, and
$g^{i}(v0)=g^{i}(v)1$, $2^{n}\le i < 2^{n+1}$.
It is clear that $g$ realises the
{\em adding machine}: the result of the application
of $g$ to a branch $(a_{1},a_{2},...)$ is
the 2-adic sum $(a_{1},a_{2},...) + (1,0,...)$.
Let $G = cl(\langle g\rangle )$;
then $G\cong {\bf Z}_{2}$, the additive group of 2-adic
integers. \parskip0pt

If $g'\in G\cap (Ker\pi_{n}\setminus Ker\pi_{n+1})$
then for every level $m>n$ there exists exactly
$2^{n}$ $g'$-cycles of length $2^{m-n}$. 
This follows from the fact that $g'$ can 
be interpreted as a 2-adic number
$c_{n+1}2^{n+1}+c_{n+2}2^{n+2}+ ...$ with $c_{n+1} =1$.
(Thus $c_{n+1}2^{n+i+1}+c_{n+2}2^{n+i+2}+ ...$ 
corresponds to $(g')^{2^i}$.) \parskip0pt  

By the description of conjugacy in $Iso(T)$ given
in \cite{GNS} we see that $g_{1}$ and $g_{2}$ from $G$ have
the same isometry type if and only if there is a number
$n\in \omega$ such that
$g_{1},g_{2} \in Ker\pi_{n}\setminus Ker\pi_{n+1}$.
This obviously implies that the group $G$ has a small
number of isometry types. $\Box$

\bigskip 

\begin{thm} \label{nilp} 
Let $T$ be a locally finite rooted tree.  

(2) If a pronilpotent group $G$ has a faithful representation in $Iso(T)$ 
with a small number of isometry types, then the number of primes dividing 
$|G|$ is finite. 
\parskip0pt 

(2) If a nilpotent group $G$ has a faithful representation in $Iso(T)$ 
with a small number of isometry types, then the periodic part of $G$ 
is of finite exponent. 
\end{thm} 

{\em Proof.} (1) Let $G$ be the inverse limit of the sequence
$G_{1} \leftarrow G_{2} \leftarrow G_{3} \leftarrow...$
corresponding to an appropriate representation on a tree.
Let $P$ be the set of primes dividing $|G|$. 
Since $G$ is pronilpotent, for every
$X\subseteq P$ there exists $g_{X}\in G$ such that a prime
number $p$ divides $|g_{X}|$ if and only if $p\in X$.
If $P$ is infinite there is continuum many elements of $G$
having pairwise distinct orders in $Iso(T)$. 
Any two such elements cannot be conjugated in $Iso(T)$. 
\parskip0pt

(2) Assume the contrary. 
Taking a subgroup if necessary we may assume that $G$ is nilpotent 
and periodic. 
Note that if each factor of the upper central series is of 
finite exponent then $G$ is of finite exponent. 
Thus we can find a characteristic subgroup $G^{-} <G$ of finite exponent 
such that the center $C(G/G^{-} )$ is infinite and is not of bounded exponent. 

Let $\pi_n$ be the standard projection $G\rightarrow G_n$ corresponding 
to the restriction to $T_n$. 
We now build a sequence of natural numbers $n_1 < n_2 <...<n_i ,...$ 
and a sequence $g_1 ,..., g_i ,... \in G\setminus G^{-}$ 
with $g_i G^{-} \in C(G/G^{-})$ such that for all 
$i\in\omega \setminus \{ 0\}$ the following properties hold: 
\begin{quote}  
(a) $\prod_{j\le i}|g_j  | < |g_{i+1} |$ ; \\ 
(b) $|\pi_{n_i} (g_i )| = |g_i |$ ; \\ 
(c) $g_{i+1} \in Ker (\pi_{n_i})$.  
\end{quote} 
This can be done by induction as follows. 
Using the fact that $C(G/G^{-})$ is infinite and is not of finite 
exponent we choose $g_{i+1}$ such that conditions (a) and (c) 
are satisfied. 
Since $|g_{i+1}|<\infty$, we can find $n_{i+1}$ such that (b) 
is satisfied. \parskip0pt 

Now for every subset $Y\subset \omega$ take 
$g_Y := \prod \{ g_i : i\in Y\}$.  
It is clear that $g_Y$ is the limit of all products 
$\prod \{ g_j : j\le i, j\in Y\}$. 
If $Y\not= Y'$ and $i$ is the least nutural number 
from $Y\Delta Y'$, then the restrictions of $g_Y$ and $g_{Y'}$ 
to $T_{n_i}$ have distinct orders and cannot be conjugated 
by an isometry. 
This contradicts the assumption that $G$ has a small number of isometry types.   
$\Box$ 
\bigskip 

Note that it is not very difficult to realize the situation 
of Theorem \ref{nilp}(2). 
The following observation gives a wide class of examples 
of representations with small numbers of isometry types.

\begin{prop} \label{product} 
Let $H$ be a finite group. 
Then the product $G= H^{\aleph_{0}}$ has a faithful action 
on a rooted tree with countably many isometry types.
\end{prop}

{\em Proof.} Consider the inverse system of projections
$$
H \leftarrow H\oplus H \leftarrow H\oplus H\oplus H
\leftarrow...$$
where the groups consist of sequences
$$
H=\{ h_{1},h_{2},...,h_{n}\} \mbox{ , }
H\oplus H=\{ h_{1}h_{1},h_{1}h_{2},...,h_{n-1}h_{n},h_{n}h_{n}\}
$$
$$,...,H\oplus ...\oplus H= \{h_{1}...h_{1},...,h_{n}...h_{n}\},... \mbox{ . }
$$ 
The same sequences form an $n$-ary tree where the groups
act by multiplication: $h_{i}w\cdot h_{l}=h_{j}w$ for
$h_{i}\cdot h_{l}= h_{j}$,
$h_{i}h_{m}w\cdot h_{k}h_{l}=h_{j}h_{n}w$ for
$(h_{i},h_{m})\cdot (h_{k},h_{l})= (h_{j},h_{n})$,... .
It is easy to see that if $g\in G$ and $|g|=s$ then there
are sequences $n_{1},n_{2},...,n_{t}$ and $s_{1},...,s_{t}=s$
such that $|\pi_{n_{1}-1}(g)|=1$,
$|\pi_{n_{1}}(g)|=s_{1}=|\pi_{n_{2}-1}(g)|$,
$|\pi_{n_{2}}(g)|=s_{2}=|\pi_{n_{3}-1}(g)|,...$ and
$|\pi_{n_{t}}(g)|=s_{t}= |\pi_{m}(g)|$ for all $m\ge n_{t}$.
In that case for every level $n_{j}\le m< n_{j+1}$
all $g$-cycles are of length $s_{j}$
and for every level $m\ge n_{t}$ all $g$-cycles are of length
$s$. By the description of conjugacy in $Iso(T)$ given in
[GNS] we see that $g_{1}$ and $g_{2}$ from $G$ have
the same isometry type if and only if the
corresponding sequences for them are the same.
This obviously implies that the group $G$ has a countable
number of isometry types. $\Box$

\bigskip
We now see that any group of the form
$$
\mathbb{Z}(n_{1})^{\omega}\oplus...\oplus \mathbb{Z}(n_{k})^{\omega}
\oplus \mathbb{Z}_{p_{1}}^{k_1}\oplus...
\oplus \mathbb{Z}_{p_{l}}^{k_l}
$$
has a faithful representation with a small number of isometry
types. 

{\bf Conjecture.} An abelian profinite group $A$ has a faithful 
representation with a small number of isometry types if and only if
$$
A\cong \mathbb{Z}(n_{1})^{\omega}\oplus...\oplus
\mathbb{Z}(n_{k})^{\omega} \oplus \mathbb{Z}_{p_{1}}^{k_1}\oplus...
\oplus \mathbb{Z}_{p_{l}}^{k_l}.
$$

\bigskip

Proposition \ref{proptrans} motivates the question of description 
of subgroups $G$ of $Iso (T)$ having a small number of isometry types 
in $G\times G$. 
We have a small remark concerning this case. 

\begin{lem} 
If $G<Iso(T)$ has a small number of isometry types of $G\times G$, 
then $G$ does not have an element of infinite order. 
\end{lem}

{\em Proof.} 
Indeed, if $|g|=\infty$, then for any $z_1 ,z_2 \in \mathbb{Z}_p$ 
with $z_1 \not=z_2$ we have $g^{z_1} \not=g^{z_2}$. 
If an isometry $h\in Iso(T)$ centalizers $g$, then it must centralize 
any $g^{z}$, $z\in\mathbb{Z}_p$. 
Thus all paires $(g,g^{z})$, $z\in \mathbb{Z}_p$, have pairwise distinct 
isometry types.  
$\Box$
\bigskip

{\bf Remark.} 
The group $A= \prod_{\omega}\langle g_{i}\rangle$
with $\langle g_{i}\rangle\cong \mathbb{Z}(p)$, is
an example of a small group in the sense of Newelski
(this was observed by G.Bezulski in his master thesis). 
It is worth noting that the action of $A$ on a binary tree 
as in Proposition \ref{product}, has a small number of 
isometry types. 
We thus consider the inverse system of projections
$$
\mathbb{Z}(2) \leftarrow \mathbb{Z}(2)\oplus \mathbb{Z}(2) \leftarrow
\mathbb{Z}(2)\oplus \mathbb{Z}(2)\oplus \mathbb{Z}(2) \leftarrow...
$$
where the groups consist of sequences
$0,1, 00,01,10,11, 000,...,111,...$ .  
The same sequences form a binary tree where the groups
act as follows: $0w+0=0w$, $0w+1= 1w$,
$00w+00=00w$, ..., $01w+11= 10w$, $11w+11= 00w$,... .
It is easy to see that if
$g'\in A\cap (Ker\pi_{n}\setminus Ker\pi_{n+1})$
then for every level $m>n$ there exists exactly
$2^{m-n-1}$ $g'$-cycles of length $2$. 
By the description of conjugacy in $Iso(T)$ given in \cite{GNS} 
we see that $g_{1}$ and $g_{2}$ from $A$ have
the same isometry type if and only if there is a number
$n\in \omega$ such that
$g_{1},g_{2} \in Ker\pi_{n}\setminus Ker\pi_{n+1}$.
This obviously implies that the group $A$ has a small
number of isometry types. \parskip0pt

On the other hand by Lemma \ref{8} the group $A$ does
not have a small number of isometry types under
the following presentation in $T= p^{{\bf N}}$: for any
$i$ the element $g_{i}$ stabilizers $T\setminus T_{1...10}$
and $(1...10jw)g_{i}= 1...10(j+1)w$). 
$\Box$

\bigskip

\subsection{Isometry groups having a small number of isometry 
types for all isometry representations} 

\begin{lem} \label{lemnilp} 
(1) Let a profinite group $G$ and a sequence $g_i\in G$, 
$i\in\omega$, satisfy the condition that for every 
$i\in\omega$ there is a continuous epimorphism 
$\rho_i$ from $G$ to a finite group such that 
$\rho_i (g_i )\not= 1$, but all $g_j$ with $j\not= i$ 
belong to the kernel of $\rho_i$.  
Then there is a locally finite rooted tree $T$ such that 
$G$ has a faithful representation in $Iso (T)$ 
with continuum many isometry types. \parskip0pt 

(2) If $G$ is a pro-$p$-group as above, then there is 
a faithful representation of $G$ in $Iso (T)$ with continuum many 
isometry types, where $T$ is the rooted tree of valency $p+1$. 
\end{lem}

{\em Proof.} 
(1) We define $T$ as a composition of two parts. 
To define the first one assume that $G$ can be presented as 
an inverse limit of a system 
$G_0\leftarrow G_1 \leftarrow ...\leftarrow G_i \leftarrow ...$
of finite groups. 
Let $N_i$, $i\in \omega$, be the kernels of the corresponding 
canonical homomorphisms $G\rightarrow G_i$. 
Then $G$ acts by left multiplication on the tree $T_C (G)$ of 
left cosets of all $N_i$, $i\in\omega$. \parskip0pt  

The second part $T_S (G)$ of the tree consists of all elements of 
$\bigcup\{ (\bigoplus^{j}_{i=0} \rho_i (G)) : j\in\omega \}$, 
where the ordering of the tree is defined by extension. 
We define an action of $G$ on this tree as follows. 
For every $g\in G$ and every vertex of the form 
$(1,...,1,\rho_{i-1}(h),\rho_{i}(h'),\rho_{i+1}(h''),...)$ 
with $h\not=1$, define the $g$-image of this vertex as 
$(1,...,1,\rho_{i-1}(h),\rho_{i}(gh'),\rho_{i+1}(h''),...)$. 
\parskip0pt 

Adding to the disjoint union $T_C (G)\cup T_S (G)$ an additional 
root $\emptyset$ together with edges to the roots of these trees 
(denoted by $\emptyset_C$ and $\emptyset_S$) we obtain the required 
tree $T$.  
Since the $G$-set $T_C (G)$ is equivariantly embedded into $T$ 
the latter is a faithful $G$-set. \parskip0pt 

For any subset $J=\{ j_1 ,j_2 ,...\}\subset \{ 1,2,...\}$ 
consider the sequence $g_{j_1}\cdot g_{j_2}\cdot ...\cdot g_{j_k}$, 
$k\in\omega$. 
Since $G$ is compact, there is the limit $g_J$ of this sequence. 
By the definition of the action of $G$ on $T_S (G)$ we easily 
see that for $J\not= K$ the elements $g_{J}$ and $g_{K}$ have 
distinct isometry types.   

(2) The proof of this statement is similar to the proof above. 
We modify it as follows. 
Replace $T_C (G)$ above by any faithful representation of $G$ 
on the rooted tree $T(p)$ of valency $p+1$ 
(see Proposition 2 of \cite{Grig}). 
It is well-known that any finite $p$-group can be embedded 
into the isometry group of a finite part of $T(p)$ 
(consisting of several levels).   
Thus every $\rho_i (G)$ can be presented as an isometry 
group of some $T(p)_{k_i}$ (consisting of $k_i$ levels). 
We now replace $T_S (G)$ above by the tree $T(p)$. 
We define an action of $G$ on this tree as follows. 
For every $g\in G$ and every vertex of the form 
$(0^{k_1},...,0^{k_{i-1}},1 ,t_1,t_2 ,...,t_{k_i} ,...)$ 
define the $g$-image of this vertex as 
$(0^{k_1},...,0^{k_{i-1}},1, \rho_{i}(g)t_1 ,\rho_{i}(g)t_2 ,...,\rho_{i}(g)t_{k_i},...)$ 
(we define the action of $g$ identically on other elements). 
The rest is obvious. 
$\Box$

\begin{thm} \label{nilp2}
Let $G$ be a nilpotent profinite group. 
If all representations of $G$ in isometry groups 
of locally finite trees have small numbers of isometry 
types, then $G$ is a group of finite rank. 
\end{thm} 

{\em Proof.}
Assume that the rank of $G$ is infinite. 
By Theorem \ref{nilp} (1) there is a prime number $p$ 
such that the Sylow $p$-subgroup of $G$ has infinite rank. 
Assume that $G$ is a pro-$p$-group. 
By composing trees as in the proof of Lemma \ref{lemnilp} we can 
reduce the situation to this case.  
\parskip0pt 

Take the upper central series $1=G_0 <G_1 <...<G_n =G$. 
It is clear that all $G_i$ are closed. 
Since $rk(G)=\infty$ we  have some $i$ with infinite $rk(G_{i+1}/G_i )$. 
Find the least $i\le n$ such that all $rk(G_{i+j+1}/G_{i+j} )$ are finite. 
We denote this number by $m$. 
Thus $rk(G/G_m )$ is finite and the group $G_m /G_{m-1}$ 
is a $\mathbb{Z}_p$-module which does not have a finite generating set.  
\parskip0pt 

Let cosets $g'_1 G_m,...,g'_t G_m$ generate topologically $G/G_m$ 
and $K$ be the intersection of $G_m$ with the closure of 
the subgroup $\langle g'_1 ,...,g'_t\rangle$.  
Note that the closure of $\langle g'_1 ,...,g'_t \rangle G_{m-1}$ 
is a normal subgroup of $G$ such that the corresponding quotient 
is isomorphic to $G_m /KG_{m-1}$ 
($=cl(\langle g'_1 ,...,g'_t\rangle )G_m /cl(\langle g'_1 ,...,g't\rangle )G_{m-1}$). 
Thus to prove the statement it suffices to build a sequence 
$g_1 ,g_2 ,..., g_i ,...\in G\setminus G_{m-1}$ 
representing pairwise distinct cosets of $G_m /KG_{m-1}$ which satisfy 
the conditions of Lemma \ref{lemnilp}(1) with respect to this group. 
\parskip0pt 

Assume that $g_1 ,g_2 ,...,g_l$ and the corresponding $rho_i$ 
are already defined. 
Since $G_m /KG_{m-1}$ is a $\mathbb{Z}_p$-module which does not have 
a finite generating set, there is an open normal subgroup $L$ of $G_m$ 
containing $KG_{m-1}$ such that $L/KG_{m-1} <(\bigcap^{l}_{i=1} ker(\rho_i ))$
$\bigcap^{l}_{i=1} ker(\rho_i )/L \not\le \langle g_1 ,...,g_l \rangle L$.  
Let $g_{l+1}$ represent an $L$-coset of $\bigcap^{l}_{i=1} ker(\rho_i )/L$ 
which does not belong to $\langle g_1 ,...,g_l \rangle L$ and let 
$\rho_{l+1}$ be the natural homomorphism onto 
$$
(\bigcap^{l}_{i=1} ker(\rho_i ))/(\bigcap^{l}_{i=1} ker(\rho_i )\cap \langle g_1 ,...,g_l \rangle L).
$$  
It is easy to see that this construction produces a required sequence. 
$\Box$ 

\bigskip

\begin{thm} \label{Z_p}
(1) Let $\Delta$ be a set of prime numbers. 
For each $p_i\in\Delta$ let $G_i$ be a pro-$p_i$-group such that all 
faithful representations of $G_i$ in isometry groups of locally finite 
rooted trees have a small number of isometry types.   
The set $\Delta$ is finite exactly when every faithful representation 
of $\oplus_{p_i\in\Delta}G_{i}$ in the isometry group of a locally finite 
rooted tree has a small number of isometry types. 

(2) An abelian profinite group $G$ has a small number of isometry types 
for all faithful representations of $G$ in isometry groups of locally 
finite rooted trees if and only if the set of prime numbers dividing 
$|G|$ is finite and $G$ is of finite rank. 
\end{thm}

{\em Proof.} (1) Necessity of the statement follows from Theorem \ref{nilp}. 
Consider the case of a representation of $H \oplus G_s$, 
where $H$ is a pro-$\Delta$-group where $\Delta$ is finite and $p_s \not\in\Delta$.  
Assume that all representations of $H$ as closed subgroups of isometry groups 
of locally finite rooted trees have a small number of isometry types. 
We want to prove that so does the group $H \oplus G_s$. 
\parskip0pt 

Let $H \oplus G_s$ act on $T$ as a closed subgroup of $Iso(T)$.  
Let $g \in G_s$.  
Consider the tree of $g$-cycles $T/\langle g\rangle$ 
induced by $g$ on $T$. 
To each vertex $\langle g\rangle t$ of $T/\langle g \rangle$ 
we assign two numbers: $v(t)$, the valency of $t$, and 
$m(t) =|\langle g \rangle t|$. 
Since $G_s$ has a small number of isometry types, 
there is a small number of such trees corresponding to elements 
of $G_s$ (up to label preserving isomorphism). 
Note that $m(t)$ is a power of $p_s$. \parskip0pt 
   
Since $[H ,G_s ]=1$ the action of $H$ on $T$ induces 
an action of $H$ on $T/\langle g\rangle$. 
If $h\in H$ defines a $h$-cycle $\langle h \rangle t$ in $T$, 
then let $l_{h}(t )$ be the length of the corresponding 
$h$-cycle in $T/\langle g\rangle$. 
Since $p_s$ does not divide $|\langle h\rangle t|$, we 
see that $l_h (t)=|\langle h\rangle t|$. 
Moreover we also see that the length of the $g h$-cycle of 
$t$ is equal to $m(t) l_h (t)$.  
Now it is straightforward that the isometry type of $h$ 
on $T/\langle g \rangle$ determines the isometry type 
of $g h$ on $T$ (even the isometry type of the pair $(g ,h)$). 
Since the number of isometry types of elements of $H$ on 
$T/\langle g\rangle$ is less than continuum we easily 
obtain the statement of the proposition. 

(2) Let $\Delta$ be a set of prime numbers. 
Let $G= \prod_{p\in\Delta}\mathbb{Z}^{l_p}_{p}$ where all $l_p$ are finite. 
It is enough to prove that if $\Delta$ is finite then 
every faithful representation of $G$ in $Iso(T)$ has 
a small number of isometry types. \parskip0pt

By part (1) it suffices to consider the case when $\Delta$ is a singleton. 
Let $l_p =l\ge 1$. 
Let $T$ be a tree as in the formulation and $\mathbb{Z}^{l}_{p}$ 
act as a closed subgroup of $Iso(T)$.  
By $\pi_n$ we denote the canonical homomorphism from $\mathbb{Z}^{l}_{p}$ 
to the group of all restrictions of $\mathbb{Z}^{l}_{p}$ to $T_n$, $n\in \omega$. 
Thus for every $g\in \mathbb{Z}^{l}_{p}$ every $\pi_n$ 
maps the group $g^{\mathbb{Z}_p}$ to a cyclic $p$-group. 
\parskip0pt 

We fix a tuple of free topological generators of $\mathbb{Z}^{l}_{p}$: 
$g_1 ,...,g_l$ (so that 
$\mathbb{Z}^{l}_{p}=g^{\mathbb{Z}_{p}}_{1}\times ...\times g^{\mathbb{Z}_{p}}_{l}$).  
To an element $g^{z_1}_1 \cdot ...\cdot g^{z_l}_{l}$ we 
associate the tuple $(v_p (z_1 ) ,...,v_p (z_l ))$ of the corresponding 
values with respect to the $p$-adic valuation. 
We claim that given an action of $\mathbb{Z}^{l}_{p}$ on $T$ 
{\em the tuple $(v_p (z_1 ) ,...,v_p (z_l ) )$ determines 
the size of the $g^{z_1}_1 \cdot ...\cdot g^{z_l}_{l}$-cycle 
of any $t\in T$.} 
In particular this tuple determines the isometry type of 
$g^{z_1}_1 \cdot ...\cdot g^{z_l}_{l}$ 
(thus the action has a small number of isometry types). 
We prove this by induction by $l$. \parskip0pt 

Let $l=1$. 
For any $t\in T_n$ the $g^{z_1}_1$-cycle of $t$ is determined by 
$\pi_n (g^{z_1}_{1})$. 
Thus to find $|g^{z_1}_1 t|$ we may assume that $z_1$ belongs 
to $\mathbb{Z}$.  
Let $v_p (z_1 ) = s$, i.e. $z_1 = p^{s}k$ with $(k,p)=1$.  
Since the $g^{k}_{1}$-cycle of $t$ coincides with the $g_1$-cycle 
of $t$, the length of the $g^{z_1}_1$-cycle of $t$ 
equals the first integer not less than $|\langle g_1 \rangle t|/ p^{s}$. 
The rest of this case is an obvious application of \cite{GNS}. 
\parskip0pt 

Now consider the case of an action of $H\times g^{\mathbb{Z}_p}_l$, 
where $H= g^{\mathbb{Z}_p}_1 \times ...\times g^{\mathbb{Z}_p}_{l-1}$.  
Assume that for any isometric action of $H$ on a locally finite 
rooted tree $T'$ and any $t'\in T$ the number 
$|\langle h \rangle t'|$, $h\in H$, is determined by 
the corresponding tuple of $p$-adic values of $h$. 
In particular the isometry type of $h$ is 
determined by this tuple of $p$-adic values. 
We want to prove these statement for the action on $T$ and 
any $g^{z_1}_1 \cdot ...\cdot g^{z_l}_l$. 
\parskip0pt 

As above we assume that $t\in T_n$ and $z_l =p^{s}k$ with $(k,p)=1$.  
Then $g_l$ and $g^{k}_l$ have the same cycles of $t$. 
Let $h= g^{z_1}_1 \cdot ...\cdot g^{z_{l-1}}_{l-1}$.
Consider the tree of $h$-cycles $T/\langle h\rangle$ 
induced by $h$ on $T$. 
To each vertex $\langle h\rangle t'$ of $T/\langle h \rangle$ 
we assign two numbers: $v(t')$, the valency of $t'$, and 
$m(t') =|\langle h \rangle t'|$.  
Note that $m(t')$ is a power of $p$ and is determined by 
the corresponding tuple of $p$-adic values of $h$. \parskip0pt 

Since $[H ,g_l ]=1$ the action of $g_l$ on $T$ induces 
an action of $g_l$ on $T/\langle h\rangle$. 
Then $g_l$ and $g^{k}_l$ have the same cycles of $\langle h\rangle t$. 
\parskip0pt 

As a result we see that $|\langle g^{z_l}_l\rangle t|$ equals 
the first integer not less than $|\langle g_l\rangle t|/p^s$, and 
the length of the $g^{z_l}_l$-cycle of $\langle h \rangle t$ 
in $T/\langle h\rangle$ equals the first integer not less than  
$|\langle g_l\rangle (\langle h\rangle t)|/p^s$. 
Dividing $|\langle g^{z_l}_l \rangle t|$ by the length 
of the $g^{z_l}_l$-cycle of $\langle h \rangle t$ we 
find $|(\langle h\rangle t) \cap (\langle g^{z_l}_l \rangle t)|$.
This number uniquely determines all equalities $h^i (t) = g^{jz_l}_l (t)$ 
for $i\le |\langle h\rangle t|$, $j\le |\langle g^{z_l}_l \rangle t|$. 
In particular it uniquely determines the length of 
the $hg^{z_l}_l$-cycle of $t$. 
$\Box$

\bigskip 

\begin{quote} 
{\bf Conjecture.} 
Let $G$ be a closed subgroup of $GL_d (\mathbb{Z}_p )$. 
Then $G$ has a small number of isometry types for all faithful 
actions on locally finite rooted trees. 
\end{quote} 
In the case of abelian groups the conjecture holds by 
Theorem \ref{Z_p} (2). 
By Lemma \ref{lemnilp} we see that groups of the form 
$H^{\omega}$ always have a representation with continuum 
of isometry types. 
On the other hand it is easy to see that this lemma 
cannot be applied to closed subgroups of $GL_d (\mathbb{Z}_p )$. 
Thus our conjecture looks reasonable. 
Another interesting question is to prove 
the converse of Theorem \ref{nilp2}.

\section{Rooted groups interpretable in $\mathbb{Q}_p$} 

One of the main obstacles to deeper study 
of isometry groups of rooted trees is the fact that 
there is no interesting first-order structure 
where these groups can be interpreted. 
To remedy this situation one can choose a rich 
algebraic object and consider isometry groups 
interpretable in it. 
Then the general case can be somehow approximated 
by information obtained in this way. 
We now describe one of possible realizations of this idea. 
\parskip0pt 
 
Let $T$ be a locally finite rooted tree. 
Let $\Gamma$ be a closed subgroup of $Iso(T)$. 
We say that $g,g'\in Iso(T)$ have the same 
$\Gamma$-type if the exists $h\in \Gamma$ 
such that $g^{h}=g'$.   
It makes sense to study $\Gamma$-types of $Iso(T)$ 
assuming that $\Gamma$ is interpetable in some natural 
model-theoretic objects.  
In this situation we say that $G\le Iso(T)$ has a small 
number of $\Gamma$-types if $G$ meets $\le \aleph_0$ 
$\Gamma$-orbits of the conjugacy action of $\Gamma$ 
on $Iso(T)$. \parskip0pt 

In this section we concentrate on trees $T$ 
and groups $\Gamma \le Iso(T)$ which are interpretable 
in the theory of the $p$-adics $\mathbb{Q}_p$. 
This provides new tools and, moreover, gives a possibility 
of extention of the matter of the previous sections by 
some additional questions. \parskip0pt 
   
For every natural $n$ we introduce the following 
equivalence relation on $Iso(T)$: 
for $g,g'\in Iso(T)$ we define $g \equiv^{\Gamma}_{n} g'$ 
if the restrictions of $g$ and $g'$ to $T_n$ are conjugated by 
an element of the projection of $\Gamma$ to $T_n$.  
For a subgroup $G\le Iso(T)$ we denote by $c^{G}_n$ 
the number of all $\equiv^{\Gamma}_{n}$-classes meeting $G$. 
We will also study the question when the zeta function 
$\Sigma_{n\ge 0} c^{G}_{n} t^n$ is rational. 
Applying the main results of \cite{hrma} we give several 
positive examples. 
\bigskip  

We start with the group $GL_d (\mathbb{Z}_p )$ and 
the corresponding tree $T_{ad}$, i.e. the tree of  
cosets in $(\mathbb{Z}_p )^d$ with respect to all subgroups 
of the form $p^k (\mathbb{Z}_p )^d$. \parskip0pt 
 
Let $\Gamma$ be a closed subgroup of $GL_{d}(\mathbb{Z}_p )$. 
Then the group $D=(\mathbb{Z}_p )^d$ becomes 
a continuous $\Gamma \mathbb{Z}_p$-module. 
Thus $\Gamma$ acts by isometries on $T_{ad}$. 
We now give some examples which show that 
typical subgroups of $GL_d (\mathbb{Z}_p )$ 
have continuum $\Gamma$-types with respect to 
their standard action on $T_{ad}$. 

{\bf Example 1.} 
Let $\Gamma =GL_d (\mathbb{Z}_p )=G$ with $d>1$.  
Since $|\mathbb{Z}_p |=2^{\aleph_0}$, 
one can produce continuum many Jordan matrices 
over $\mathbb{Z}_p$. 

Assuming that $d=2$ and $G$ is the group $SD_2 (\mathbb{Z}_p )$ 
of all diagonal matrices over $\mathbb{Z}_p$ 
with determinant 1 we obtain an abelian 
group of isometries of the tree $T_{ad}$ 
having continuum $\Gamma$-types where 
$\Gamma =GL_d (\mathbb{Z}_p )$.  

On the other hand extending this idea we can 
prove the following statement. 

\begin{prop} \label{titsalt}
Let $\Gamma =GL_d (\mathbb{Z}_p )$ and $G\le GL_d (\mathbb{Z}_p )$. 
If $G$ is a closed subgroup of $Iso(T_{ad})$ which 
does not have a soluble subgroup of finite 
index, then $G$ meets $2^{\aleph_0}$ 
$\Gamma$-types. 
\end{prop}

{\em Proof.} 
By Tits alternative $G$ has a free subgroup. 
The proof of \cite{tits} shows that $G$ contains 
a semisimple element $g$ of infinite order. 
Then any pair of powers $g^{\mathbb{Z}_p}$ cannot 
be conjugated in $\Gamma$. 
$\Box$ 
\bigskip 

{\bf Example 2.}
Let $\Gamma =SL_2 (\mathbb{Z}_p )$ and $G=UT_2 (\mathbb{Z}_p )$ 
be the abelian group of all unitriangular matrices over 
$\mathbb{Z}_p$, $p>2$.   
Straightforward computations show that matrices 
$(a_{11},a_{12},a_{21},a_{22})$ and $(b_{11},b_{12},b_{21},b_{22})$ 
frome $G$ (written as vectors, i.e. $a_{21}=b_{21}=0$ and 
$a_{11}=b_{11}=a_{22}=b_{22}=1$) are 
conjugated by $(x_{11},x_{12},x_{21},x_{22})\in \Gamma$ if 
and only if $x_{21}=0$ and $b_{12}=x^2_{11}a_{12}$.  

It is worth noting that since  
$(\mathbb{Z}^{\times}_p )^2 =\mathbb{Z}^{\times}_p$, 
any two elements $\bar{a},\bar{b}\in G$ as above with 
$\nu (a_{12}) =\nu(b_{12})$,  
are $\Gamma$-conjugated (by a diagonal matrix). 
In this case $G$ has a small (i.e.countable) number 
of $\Gamma$-types. 

\bigskip 

{\bf Remark.} 
It is natural to ask if the actions above 
have a small number of isometry types 
with respect to the corresponding tree $T$ 
and $Iso(T)$. 
It is clear that if for some $\Gamma <Iso(T)$, 
$G$ has a small number of $\Gamma$-types 
then $G$ has a small number of isometry types. 
Thus the case when $G=UT_2 (\mathbb{Z}_p )$ acts 
on the corresponding tree $T_{ad}$ is clear.

Assume that $G$ is the group $SD_2 (\mathbb{Z}_p )$ 
of all diagonal matrices over $\mathbb{Z}_p$ 
with determinant 1. 
Then 
\begin{quote} 
{\em the natural isometric action of $G$ on 
the tree $T_{ad}$ has a small number of isometry types.}   
\end{quote} 
Indeed identify $G$ with the set of diagonals  
$\{ (u,u^{-1}): u\in \mathbb{Z}^{\times}_p \}$. 
Then $G$ acts on the module $\mathbb{Z}_p \times\mathbb{Z}_p$ 
by the action: 
$$
(u,u^{-1})\cdot (z_1 ,z_2 ) = (uz_1 ,u^{-1}z_2 ). 
$$ 
This action naturally induces the action of $G$ on $T_{ad}$. 
It is clear that $G\cong \mathbb{Z}^{\times}_p$. 
On the other hand it is known that the group $\mathbb{Z}^{\times}_p$ 
can be decomposed into a direct sum isomorphic to 
$\mathbb{Z}(p) \oplus \mathbb{Z}_p$, where the factor isomorphic 
to $\mathbb{Z}(p)$ is represented by the cyclic subgroup 
$\{ \varepsilon_1 ,...,\varepsilon_{p-1}\}$ of $p-1$-th roots of unity. 
The factor $\mathbb{Z}_p$ is represented by all elements of 
$\mathbb{Z}^{\times}_p$ which are congruent to 1 mod $p$ 
(see Chapter 18 in \cite{fuchs}). 
It only remains to apply Theorem \ref{Z_p}. 
\parskip0pt 

In fact the latter factor contains $l$ such that all other elements 
are represented by $l^{z}$, where $z\in p\mathbb{Z}_p$. 
To find this $l$ let $t$ be a primitive root mod $p$ 
which is a primitive root modulo every power of $p$ 
and let $l= t^{p-1}$ (see \cite{fuchs}, p.316). 

As a result the action of $(u,u^{-1})=(l^{z},l^{-z})$ 
on a vertex $(z_1 ,z_2 )+p^{k}\mathbb{Z}_p$ 
is defined as follows. 
Assuming that $z_1 ,z_2 \in \omega$ and $z_1 ,z_2 <p^k$ 
we represent 
$(z_1 ,z_2 )=(p^{s_1}\varepsilon_i l^{t_1}, p^{s_2}\varepsilon_j l^{t_2})$ 
with $s_1 ,s_2 <k$ and $t_1 < k-s_1 -1$ and $t_2 <k-s_2 -1$ 
(this is possible because $l^i \equiv 1$ mod $p^k$ 
if and only if $p^{k-1}|i$). 
Then the action $(u,u^{-1})\cdot ((z_1 ,z_2 )+p^{k}\mathbb{Z}_p )$ 
is defined by the representative 
$(p^{s_1}\varepsilon_i l^{z+t_1}, p^{s_2}\varepsilon_j l^{-z+t_2})$.

\bigskip

According to Proposition \ref{titsalt}, when $\Gamma$ 
and $T_{ad}$ are as in this proposition, groups with 
a small number of $\Gamma$-types must be almost soluble. 
According to Theorem 1 of \cite{Grig2} the Hilbert-Poincar\'{e} 
series of such a group $G<GL_d (\mathbb{Z}_p )$ 
is $\mathbb{Q}$-rational. 

To show that zeta-functions of the numbers of conjugacy 
classes are rational for these examples we now formulate 
some general theorem. 

We consider the $p$-adics $\mathbb{Q}_p$ with the standard 
valuation $\nu :\mathbb{Q}^{*}_{p}\rightarrow \mathbb{Z}$ 
(and the value group), 
the valuation ring $\mathcal{O}=\mathbb{Z}_p$, the maximal 
ideal $M=\mathbb{Z}_p \setminus \mathbb{Z}^{\times}_p$ 
and the corresponding residue field. 
We remind the reader that $\nu$ is a homomorphism such that 
$\nu (x+y)\ge inf (\nu (x),\nu (y))$ and $\nu (p)=1$. 
Then $\mathbb{Z}_p =\{ x\in \mathbb{Q}_p : \nu (x)\ge 0\}$ is 
the valuation ring of $\mathbb{Q}_p$ and $\mathbb{Z}/p\mathbb{Z}$ 
is the corresponding residue field. \parskip0pt 

For each $n\in \omega$ we add the sort $GL_n (\mathbb{Q}_p )$ 
and the set of lattices 
$S_n (\mathbb{Q}_p )=GL_n (\mathbb{Q}_p )/GL_n (\mathbb{Z}_p ) $. 
We also add the natural map 
$GL_n (\mathbb{Q}_p )\rightarrow S_n (\mathbb{Q}_p )$. 
Structures of this form were introduced in \cite{hahrma} 
and \cite{hrma}, where sorts of these structures were 
called geometric. 
It is proved in \cite{hrma} that $\mathbb{Q}_p$ in 
this language admits elimination of imaginaries. 
The following theorem from \cite{hrma} (Theorem 6.2) 
plays the central role in this section. 

Let $R=(R_l )_{l\in \omega^{r}}$ be a definable family of 
subsets of $\mathbb{Q}^{N}_p$ (i.e. a definable subset of 
$\mathbb{Q}^{N}_p \times \mathbb{Z}^r$, where $\mathbb{Z}$ is 
the value group and $\omega \subset \mathbb{Z}$). 
Let $E=(E_l )_{l\in\omega^{r}}$ be a definable family of 
equivalence relations on $R$ (i.e. every equivalence class 
is contained in an $l$-fibre of $R$ for some $l\in\omega^{r}$). 
Suppose that for each $l\in\omega^{r}$ the set of equivalence 
classes $R_l /E_l$ is finite. 
Let $a_l =|R_l /E_l |$. 
Then the power series 
$\sum_{l\in\omega^{r}}t^l\in \mathbb{Q}[[t_1 ,...,t_r ]]$ 
is $\mathbb{Q}$-rational. 

We now describe how it can be applied in our situation. 

\begin{prop} \label{interpretation}
Let $\Gamma$ and $G$ be subgroups of $GL_d (\mathbb{Z}_p )$ 
definable in $\mathbb{Q}_p$. 
Let $c^{G}_n$ be the number of all $\equiv^{\Gamma}_{n}$-classes 
meeting $G$ defined with respect to $T_{ad}$. 
Then the zeta function $\Sigma_{n\ge 0} c^{G}_{n} t^n$ 
is $\mathbb{Q}$-rational.
\end{prop}

{\em Proof.} 
We interpret $T_{ad}$ in $\mathbb{Q}_p$ as the set 
of equivalence classes on $\mathbb{Z}^d_p \times \omega$  
(where $\omega \subset \mathbb{Z}$ is the non-negative part 
of the value group) with respect to the following equivalence relation. 
Let 
\begin{quote} 
$(\bar{z}_1 ,c_1 )\sim (\bar{z}_2 ,c_2 )$ if 
$(c_1 = c_2 )\wedge (\nu (\bar{z}_1 -\bar{z}_2 )\ge c_1 )$. 
\end{quote} 
It is clear that this defines the set of cosets of all groups 
of the form $p^c (\mathbb{Z}_p )^d$. 
The inclusion relation of cosets (which is the order relation 
of $T_{ad}$) is defined by representatives: 
$$
( \bar{z}_1 ,c_1 )\le (\bar{z}_2 ,c_2 ) \Leftrightarrow  
(c_1 \ge c_2 )\wedge (\nu (\bar{z}_1 -\bar{z}_2 )\ge c_1 ). 
$$
Let $R= (R_i )_{i\in\omega} = G\times \omega$. 
We define a family $E=(E_i )$ of equivalence relations 
on $R$ as follows: 
$$ 
((g_1 ,c_1 ),(g_2 ,c_2 ))\in E \Leftrightarrow 
(c_1 =c_2 )\wedge (\exists \gamma \in \Gamma ) 
$$
$$ 
(\forall \bar{z}\in (\mathbb{Z}_p )^d )(\forall c \le c_1 ) 
(\gamma (g_1 ([(\bar{z},c)]_{\sim} ))= g_2 (\gamma ([(\bar{z},c)]_{\sim} )). 
$$
The rest follows by Theorem 6.2 of \cite{hrma}. 
$\Box$
\bigskip

\begin{prop} 
Let $G$ be a finitely generated, torsion free, nilpotent 
pro-$p$-group. 
Then for some natural $d$ there is an embedding of $G$ into 
$GL_{d}(\mathbb{Z}_p )$ such that the corresponding action 
$(G,T_{ad})$ is definable in $\mathbb{Q}_p$.  
The zeta function $\Sigma_{n\ge 0} c^{G}_{n} t^n$ 
of this action with respect to $\Gamma =GL_d (\mathbb{Z}_p )$ 
is $\mathbb{Q}$-rational.
\end{prop} 

{\em Proof.} 
Let $a_1 ,...,a_n$ be a Malcev basis of $G$ (i.e. any 
element of $G$ can be written uniquely in the form 
$a^{\lambda_1}_1 \cdot ...\cdot a^{\lambda_n}_n$, 
$\lambda_i \in \mathbb{Z}_p$). 
A well-known theorem (see \cite{merzlakov}) states 
that group multiplication and inversion in $G$ are 
given by polynomials in $\lambda_i$ with coefficients 
in $\mathbb{Q}$ as is the map 
$G\times \mathbb{Z}_p \rightarrow G$, 
$(g,\lambda )\rightarrow g^{\lambda}$. 
This gives an interpretation of $G$ in $\mathbb{Q}_p$ as 
$\mathbb{Z}^{n}_p$. \parskip0pt 

By Theorem 59.2.1 from \cite{merzlakov} there is a natural 
number $d$ and a polynomial map 
$\phi :\mathbb{Z}^{n}_p \rightarrow UT_d (\mathbb{Z}_p )$ 
giving an isomorphic representation of $G$ on $\mathbb{Q}_p$. 
Moreover the converse mape $\phi^{-1}$ is linear. 
This clearly provides a definable action $(G,T_{ad})$. 
By Proposition \ref{interpretation} we have 
the last statement of our proposition. 
$\Box$

\bigskip 

One of the main examples of $G$-trees in geometric group 
theory is the construction of $SL_2 (\mathbb{Q}_p )$-tree of lattices 
of Serre and Tits.  
We now show that the actions studied in the proposition above 
(and in Proposition \ref{interpretation}) arise as point-stabilizers 
in this example. 
The following well-known construction is taken from \cite{Serre}). 

A {\em lattice} in $\mathbb{Q}_p \times\mathbb{Q}_p$ 
is a two-generated $\mathbb{Z}_p$-submodule of 
$\mathbb{Q}_p \times\mathbb{Q}_p$ which generates the 
$\mathbb{Q}_p$-vector space $\mathbb{Q}_p \times\mathbb{Q}_p$. 
By $\mathcal{X}$ we denote the set of equivalence classes 
of lattices with respect to the equivalence 
$$ 
L\sim L' \leftrightarrow (\exists z\in \mathbb{Q}^{*}_p )(L'=Lz ).
$$ 
The distance $d([L]^{\sim} ,[L']^{\sim} )$ between the corresponding 
classes is defined as follows.  
Find $L''\sim L'$ such that 
$L''=\langle e_1 p^{a}, e_2 p^{b}\rangle < \langle e_1 ,e_2 \rangle =L$. 
Then the number $|a-b|$ does not depend on the choice of the basis 
$e_1 ,e_2$ and the representative $L''$; so it is taken as the 
corresponding distance. \parskip0pt 

The metric space $(\mathcal{X},d)$ is a simplicial tree of valency $p+1$, 
and $PGL_2 (\mathbb{Q}_p )$ has an isometric action on $\mathcal{X}$. 
Moreover the group $PSL_2 (\mathbb{Q}_p )$ acts on $\mathcal{X}$ 
without inversions.  
To interpret this $G$-space in $Th(\mathbb{Q}_p )$, the set 
of lattices is usually identified with the left coset space 
$S_2 = GL_2 (\mathbb{Q}_p )/GL_2 (\mathbb{Z}_p )$ (as above, 
see \cite{hahrma}).  
Then $\mathcal{X}$ can be considered as the coset space 
$GL_2 (\mathbb{Q}_p )/(\mathbb{Q}^{*}_p\cdot GL_2 (\mathbb{Z}_p ))$,  
where $PGL_2 (\mathbb{Q}_p )$ ($PSL_2 (\mathbb{Q}_p )$ resp.)  
acts from the left and $\mathbb{Q}^{*}_p$ now stands for the 
subgroup of scalar matrices. 

Formally we consider $\mathcal{X}$ as $GL_2 (\mathbb{Q}_p )$ with 
respect to the following equivalence relation: 
$$
(a_{ij})_{i\le 2,j\le 2} \sim (b_{ij})_{i\le 2,j\le 2} \Leftrightarrow 
(\exists (c_{ij})_{i\le 2,j\le 2},(d_{ij})_{i\le 2,j\le 2} \in 
GL_2 (\mathbb{Q}_p )) (\exists q\in \mathbb{Q}^{*}_p ) 
$$
$$
[\nu (c_{ij})\ge 0 \wedge \nu (d_{ij})\ge 0 \wedge 
((c_{ij})_{i\le 2,j\le 2}\cdot (d_{ij})_{i\le 2,j\le 2} = E )\wedge  
$$
$$
((a_{ij})_{i\le 2,j\le 2} = (b_{ij})_{i\le 2,j\le 2}\cdot (q\cdot c_{ij})_{i\le 2,j\le 2})]. 
$$
The group $PGL_2 (\mathbb{Q}_p )$ ($PSL_2 (\mathbb{Q}_p$ resp.) 
is considered as $GL_2 (\mathbb{Q}_p )$ ($SL_2 (\mathbb{Q}_p$) 
with respect to the equivalence relation induced by the center. 
The action of $PGL_2 (\mathbb{Q}_p )$ on $\mathcal{X}$ is defined 
by left multiplication of the corresponding representatives. 
\parskip0pt 

To define the distance $d$ on $\mathcal{X}$ we use 
the following fomula: 
$$ 
d((a_{ij})_{i\le 2,j\le 2}, (b_{ij})_{i\le 2,j\le 2}) =k \Leftrightarrow 
(\exists t\in \mathbb{Q}^{*}_p )(\exists (c_{ij})_{i\le 2,j\le 2})\in 
GL_2 (\mathbb{Z}_p )
$$
$$
(\exists q,r \in \mathbb{Z}_p )[k= |\nu (q)-\nu (r)| \wedge 
[(qa_{11}, ra_{12}, qa_{21},ra_{22}) = (b_{ij})_{i\le 2,j\le 2}\cdot 
(t\cdot c_{ij})_{i\le 2,j\le 2}]. 
$$ 

Let $L_0 =\langle (1,0),(0,1)\rangle$. 
The stabilizer of the point  $[L_0 ]^{\sim}$ is exactly 
the subgroup of $PGL_2 (\mathbb{Q}_{p})$ consisting 
of all cosets of the form $g\mathbb{Q}^{*}_p$ where 
$g\in GL_2 (\mathbb{Z}_p)$. 

Each vertex of $\mathcal{X}$ is represented by 
a unique lattice $L\subseteq L_0$ such that 
$L_0 /L \cong \mathbb{Z}_p /p^{n}\mathbb{Z}_p$, 
where $n$ is the distance between $L_0$ and $L$. 
In this case $L/p^n L_0$ is a direct factor of rank 1 of 
the $\mathbb{Z}_p /p^n \mathbb{Z}_p$-module $L_0 /p^{n}L_0$. 
When $L = \langle u,vp^{n}\rangle$ for some $u$ and $v$ 
generating $L_0$, we identify $L$ with the 
$\mathbb{Z}_p /p^n \mathbb{Z}_p$-module generated by $u+p^n L_0$. 

Formally we present the set of these modules as follows. 
In the rooted tree $\mathcal{L}$ of all cosets of 
$L_0 = \mathbb{Z}_p \times \mathbb{Z}_p$ 
of the form $u + p^n (\mathbb{Z}_p \times \mathbb{Z}_p )$ 
with $\nu (u)=0$ we introduce the following equivalence relation :  
$$
u+p^n L_0 \approx u' +p^l L_0 \Leftrightarrow 
(n=l) \wedge (\exists x\in \mathbb{Z}^{\times}_p )(xu-u' \in p^n L_0 ). 
$$  
It is easy to see that the action of any  
$g\mathbb{Q}^{*}_p$ with $g\in GL_2 (\mathbb{Z}_p )$ on 
the set of elements of $\mathcal{X}$ which are $n$-distant 
from $L_0$, corresponds to the action of $g$ 
on the set of these classes. 
Moreover if $g$ and $g'\in GL_2 (\mathbb{Z}_p )$ determine 
the same element of $PGL_2 (\mathbb{Z}_p )$, 
then $g(u+p^n L_0 )\approx g'(u +p^n L_0 )$, where $u$ as above.  
Thus the stabilizer of the element $L_0$ in $\mathcal{X}$ 
is identified with $PGL_2 (\mathbb{Z}_p )$, which is considered 
with respect to its action on $\mathcal{L}/\approx$.  
\parskip0pt 

Note that the rooted tree $\mathcal{L}/\approx$ is divided 
into levels represented by cosets $u+p^n L_0$ with the same $p^n$. 
Moreover the vertices $u+p^n L_0$ and 
$u'+p^{n+1} L_0$ are ajacent if and only if 
$(\mathbb{Z}_p /p^n \mathbb{Z}_p )(u+p^n L_0 )=(\mathbb{Z}_p /p^n \mathbb{Z}_p )(u' +p^n L_0 )$
(by the definition of $\mathcal{X}$). 

\begin{thm}  \label{PGL}
In the notation above let $\Gamma =G= PGL_2 (\mathbb{Z}_p )$. 

(1) Let $c^{G}_n$ be the number of all $\equiv^{\Gamma}_{n}$-classes 
meeting $G$ defined with respect to the rooted tree $\mathcal{L}/\approx$. 
Then the zeta function $\Sigma_{n\ge 0} c^{G}_{n} t^n$ 
is $\mathbb{Q}$-rational.

(2) The group $G$ has a small number of isometry types with respect to 
the action on $\mathcal{L}/\approx$ and  
there is continuum $\Gamma$-types meeting $G$. 
\end{thm} 

{\em Proof.} 
(1) Repeating the proof of Proposition \ref{interpretation} 
we interpret $\mathcal{L}/\approx$ in $\mathbb{Q}_p$ as 
the set of equivalence classes on a definable subset of 
$\mathbb{Z}^2_p \times \omega$  (where $\omega \subset \mathbb{Z}$ 
is the non-negative part of the value group).  
Then we consider 
$R= (R_i )_{i\in\omega} = GL_2 (\mathbb{Z}_p )\times \omega$. 
with respect to the family $E=(E_i )$ of equivalence relations 
on $R$ as it was described above. 
The rest follows by Theorem 6.2 of \cite{hrma}. 

(2) By obvious arguments from linear algebra we may assume that 
elements of $PGL_2 (\mathbb{Z}_p )$ are presented by upper 
triangular matrices. 
Let $(a_{11}, a_{12}, a_{22})$ be a tuple of non-zero entries 
of such a matrix $g\in G$. 
Then multiplying $g^n$ with $u=(z_1 ,z_2 )\in \mathbb{Z}_p \times\mathbb{Z}_p$ 
we obtain a representative of the projective line corresponding to 
$$
(z_1 +\frac{a_{12}z_2 }{a_{11}-a_{22}}(1-(\frac{a_{22}}{a_{11}})^n), (\frac{a_{22}}{a_{11}})^n z_2 ). 
$$
Decompose the group $\mathbb{Z}^{\times}_p$ 
into a direct sum isomorphic to 
$\mathbb{Z}(p) \oplus \mathbb{Z}_p$, where the factor isomorphic 
to $\mathbb{Z}(p)$ is represented by the cyclic subgroup 
$\{ \varepsilon_1 ,...,\varepsilon_{p-1}\}$ of $p-1$-th roots of unity. 
The factor $\mathbb{Z}_p$ is represented by all elements of 
$\mathbb{Z}^{\times}_p$ which are congruent to 1 mod $p$ 
(see Chapter 18 in \cite{fuchs}). 
As above we find $l$ such that all elements of this factor 
are represented by $l^{z}$, where $z\in p\mathbb{Z}_p$ 
(see \cite{fuchs}, p.316). 
Assuming that $a_{22}/a_{11} = \varepsilon_j l^{z}$ we see that 
$(\frac{a_{22}}{a_{11}})^n \equiv 1$ mod $p^k$ if and only if 
$\varepsilon^n_j =1$ and $p^{k-1}|nz$. 

Assume $(a'_{11}, a'_{12}, a'_{22})$ be a tuple of non-zero entries 
of a matrix $g'\in G$ such that $\nu (a_{12})= \nu (a'_{12})$, 
$a'_{22}/a'_{11} = \varepsilon_j l^{z'}$ and $\nu (z) = \nu (z')$.  
Then for every pair $(z_1 ,z_2 )\in \mathbb{Z}_p \times \mathbb{Z}_p$ 
the elements $g$ and $g'$ have cycles of the same length 
at the element of $\mathcal{X}$ represented by  $(z_1 ,z_2 )+ p^m L_0$.  
This proves the statement. 
$\Box$ 

\bigskip

{\bf Remark.} There is another version of 
Proposition \ref{interpretation} which does not use 
definability in $\mathbb{Q}_p$ and which can be also 
applied to examples of this kind. 

\begin{prop} \label{interpretation2}
Let $\Gamma$ and $G$ be closed subgroups of $GL_d (\mathbb{Z}_p )$ 
. 
Let $c^{G}_n$ be the number of all $\equiv^{\Gamma}_{n}$-classes 
meeting $G$ defined with respect to $T_{ad}$. 
Then the zeta function $\Sigma_{n\ge 0} c^{G}_{n} t^n$ 
is $\mathbb{Q}$-rational.
\end{prop}

To prove this proposition we must apply a more involved argument. 
In fact we must repeat the proof of Theorem 1.2 from 
\cite{dS} with some sligth changes. 
Since this theorem is outside the main theme 
of the paper we just mention that our adaptation 
requires that formula (1.1) from \cite{dS} 
(counting the number of conjugacy classes in a finite group) 
must be replaced by the formula 
$$
c^{G}_n = |\Gamma_n |^{-1} \sum\{ |C_{\Gamma_n}(g)|:
g\in GL_d (\mathbb{Z}_p /p^n \mathbb{Z}_p ) 
$$
$$
\wedge \exists g'\in G_n (g\equiv^{\Gamma}_n g' )\} , 
$$    
where $G_n$ and $\Gamma_n$ are the natural projections 
of $G$ and $\Gamma$ with respect to 
$GL_d (\mathbb{Z}_p /p^n \mathbb{Z}_p )$. 
Since any closed subgroup of $GL_d (\mathbb{Z}_p )$ 
is compact $p$-adic analytic, the arguments of 
Section 1 of \cite{dS} show that $G$, $\Gamma$ and 
$\equiv^{\Gamma}_n$ are definable in the theory 
of analytic functions from \cite{DvdD}. 
This allows us to adapt the proof from \cite{dS} 
to our case. 

\bigskip

\section{Branch groups. Model theoretic aspects. }

As we have noted above, branch groups have a large 
number of isometry types. 
Can they be defined (in $\mathbb{Q}_p$) as closed 
subgroups of matices over $\mathbb{Z}_p$ ? 
This is not true. 
Indeed, let $G$ be the profinite completion of 
Grigorchuk's 3-generated 2-group from \cite{Grig1}. 
In fact it is shown in \cite{Grig2} that the coefficients 
$a_n$ of the Hilber-Poincar\'{e} series of $G$ do not have 
a polynomial growth. 
By Theorem 1 of \cite{Grig2} (or Interlude A of \cite{DdSMS}) 
the group $G$ is not a closed subgroups of matices over $\mathbb{Z}_p$. 
\parskip0pt 

How do branch groups look with respect to other model theoretic 
properties of profinite groups studied so far ? 
In this section we concentrate on $\omega$-categoricity 
of branch groups in the language introduced by Z.Chatzidakis 
in \cite{zoe}. 

Let $G$ be a profinite group. 
As in \cite{zoe} to $G$ we associate a structure $S(G)$ 
of the language $L=(\le ,\sim ,C^2 ,P^3 ,1)$ as follows. 
The structure $S(G)$ is defined on the set of all 
cosets $gN$ for all open normal subgroups $N<G$. 
The symbol $P$ is interpreted as follows: 
$$
S(G) \models P(g_1 N_1 ,g_2 N_2 , g_3 N_3 )\Leftrightarrow 
(N_1 =N_2 =N_3 )\wedge (g_1 g_2 N_1 = g_3 N_1 ). 
$$ 
The symbol $C$ corresponds to inclusion: 
$$
S(G)\models C(gN, hM) \Leftrightarrow (N\subseteq M )\wedge (gM=hM).
$$ 
The relation $gN\le hM$ means $N\subseteq M$ and we define 
$gN\sim hM\Leftrightarrow (gN\le hM )\wedge (hM\le gN)$.  
The costant $1$ corresponds to $G$.  

It is convenient to view $S(G)$ as an $\omega$-sorted 
structure with respect to sorts $S_n =\{ gN :|G:N|=n\}$, $n\in \omega$. 
In this language the class of structures $S(G)$ becomes elementary. 
In Section 1.5 of \cite{zoe} some natural axioms of this class 
are given: it is shown that for any structure $S$ satisfying 
these axioms the $\sim$-classes naturally form a projective system 
of finite groups such that $S$ is $S(G)$ for the projective 
limit $G$ of this system. \parskip0pt 

We say that $A\subseteq S(G)$ {\em forms a substructure} if 
$(\forall x,y\in S(G))((x\le y)\wedge (x\in A)\rightarrow y\in A)$ 
and $(\forall x,y\in A)(\exists z\in A)((z\le y)\wedge (z\le y))$. 
It is straightforward that any substructure $A\subseteq S(G)$ 
defines a profinite group $H$ with $A=S(H)$ such that 
$H$ is a continuous homomorphic image of $G$.    
On the other hand if $\phi :G\rightarrow H$ is a continuous 
epimorphism of profinite groups, then the map 
$S(H)\rightarrow S(G)$: $gN \rightarrow \phi^{-1}(gN)$ 
defines an embedding of $S(H)$ into $S(G)$ as a substructure.  
The following property can be taken as the definition of 
$\omega$-categoricity (see Section 1.6 of \cite{zoe}). 

\begin{definicja} 
We say that a countable structure $S(G)$ is $\omega$-categorical 
if for all $n, j_1 ,...,j_k \in \omega$ the group $Aut(S(G))$ 
has finitely many orbits on $(S_{j_1} \cup ...\cup S_{j_k})^{n}$ 
(i.e. $Aut(S(G))$ is {\em oligomorphic}).  
\end{definicja} 

It is easy to see that 
\begin{quote} 
for every finitely generated profinite group $G$ each sort 
$S_i (G)$ is finite (see \cite{DdSMS}, Proposition 1.6); 
i.e. $S(G)$ is $\omega$-categorical. 
\end{quote} 
Z.Chatzidakis has noticed in \cite{zoe} (Theorem 2.3) that 
when $G$ has the {\em Iwasawa Property} (IP), the structure 
$S(G)$ is $\omega$-categorical. 
We remind the reader that $G$ has (IP) if for every epimorphism 
$\theta :H\rightarrow K$ of finite groups with $H\in Im(G)$, 
and for every epimorphism $\phi :G\rightarrow K$ there is 
an epimorphism $\psi :G\rightarrow H$ such that 
$\phi =\theta \cdot \psi$.  
Note that (IP) is very close to projectivity 
(which in the case of pro-$p$-groups (IP) is equivalent to  
$p$-freeness, see \cite{FJ}, Chapter 20). 

It is worth noting that all closed subgroups of 
$GL_n (\mathbb{Z}_p )$ have $\omega$-categorical 
structures $S(G)$. 
This follows from the fact that they have a finite number 
of open subgroups of index $k$, $k\in \omega$ 
(moreover the corresponding function counting the number 
of subgroups is polynomial \cite{DdSMS}). 
Thus groups considered in the previous sections 
usually had $\omega$-categorical structures $S(G)$ 
(excluding $\mathbb{Z}(p)^{\omega}$). 

In the following proposition we consider 
{\em just infinite profinite groups}. 
This means that each proper continuous homomorphic image 
of the group is finite. 
Just infinite branch groups are described in Theorem 4 of 
\cite{Grig}. 
It states that a profinite group with a branch structure 
$\{ L_i ,H_i :i\in \omega\}$ is just infinite if and only 
if the commutator subgroup $[L_i ,L_i ]$ is of finite 
index in $L_i$ for all $i$.

\begin{prop} \label{justinf}
Let $G<Iso(T)$ be a weakly branch profinite group such that 
the restricted stabilizer $Rs(n)$ is level-transitive 
on every subtree of level number $n$. 

(1) Then $S(G)$ is $\omega$-categorical if and only if every 
sort $S_k$ is finite. 

(2) If $G$ is additionally a branch group which is just infinite,  
then $S(G)$ is $\omega$-categorical.    
\end{prop}

{\em Proof.} 
(1) The first statement of the proposition is a straightforward 
consequence of Theorem 7.5 and Proposition 8.1 of \cite{LN}. 
Indeed each automorphism $\phi$ of $S(G)$ naturally 
defines a continuous automorphism of the group $G$.  
By \cite{LN} under our assumptions, each automorphism 
of $G$ is induced by conjugation of an isometry $h$ of 
$T$ such that $G^h = G$. 
This means that if an open subgroup $N<G$ is the pointwise 
stabilizer in $G$ of the first $l$ levels of $T$, 
then $\phi (N)=N$. 
Since each open subgroup of $G$ contains such a subgroup 
$N$ we easily see that all $Aut(S(G))$-orbits on 
$S(G)$ are finite.

(2) Let us show that any open normal subgroup $N<G$ has 
a finite orbit with respect to automorphisms of $G$ 
induced by conjugation by isometries of $T$.
Let groups $L_i$ and $H_i$ be defined as in 
the definition of branch groups (see Section 1.1). 
Since $G$ is just infinite, by Theorem 4 of \cite{Grig} 
the commutator subgroup $[L_i ,L_i ]$ is of finite 
index in $L_i$, and thus $[H_i ,H_i]$ is of finite 
index in $G$. 
On the other hand for every $n$ there is a number 
$l=l(n)$ such that every subgroup of $G$ of index $n$
contains an element which does not fix the subtree $T_l$ 
pointwise. \parskip0pt 

Assume $|G:N|=n$ and $N$ is normal in $G$. 
Find an element $g\in N$ which does not fix 
the subtree $T_l$ pointwise for $l=l(n)$. 
By the proof of Theorem 4 the group $N$ contains 
the commutator subgroup $[H_{i+1} ,H_{i+1}]$. 
Since the index $|G:[H_{i+1} ,H_{i+1} ]|$ is finite 
the number of normal subgroups of $G$ of index $n$ 
is finite. 
$\Box$  

\bigskip 

{\bf Example.} 
Let $G$ be the profinite completion of 
Grigorchuk's 3-generated 2-group from \cite{Grig1}. 
Since $G$ is a finitely generated profinite group, 
$S(G)$ is $\omega$-categorical.  
It is worth noting that it is proved in \cite{Grig} that 
$G$ is a just infinite branch group. 

\bigskip 

In the remained part of the section we study 
$\omega$-categoricity of {\em wreath branch groups} 
(see \cite{LN}).  
Since it can happen that such a group is not finitely 
generated and not just infinite, Proposition \ref{justinf} 
does not work very well here.

Let $A_i$, $i\in \omega$, be the sequence of alphabets 
defining $T$ as in Introduction. 
Assume that for every $n$ we have chosen a permutation 
group $P_i \le Sym (A_i )$. 
Then the subgroup of all labellings $\{ \gamma (v)\}_{v\in T}$ 
with $\gamma (v)\in P_{i+1}$ for $v \in layer_i$ is  
{\em the iterated wreath product} of $(P_i ,A_i )$, 
$i\in \{ 1,2,...\}$ and is denoted by $\wr_{\omega} (P_i ,A_i )$.  
It is noticed in \cite{LN} (Theorem 8.2) that this group 
is a branch group. 
Thus it is called in \cite{LN} a wreath branch group. 

\begin{prop}
(1) Assume that all groups $(P_i ,A_i )$ are simple 
transitive permutation groups. 
Then the group $G= \wr_{\omega} (P_i ,A_i )$ has 
the Iwasawa property (IP) and $S(G)$ is $\omega$-categorical.  

(2) 
The group $Iso(T)$ (i.e. $\wr_{\omega} (Sym (A_i ),A_i )$ ) 
does not have the Iwasawa property (IP) and $S(G)$ is not 
$\omega$-categorical.  
\end{prop} 

{\em Proof.} 
(1) This statement follows from the observation 
that every closed normal subgroup of $G$ is of 
the form $Ker \pi_n$ for appropriate $n$. 
We think that this fact is folklor, but we 
give some scketch of it. 
Let $K$ be a closed normal subgroup of $G$. 
Let $n$ is the least natural number such that  
there is $g\in K\setminus \{ 1\}$ of {\em depth} 
$n$, i.e. $n$ is the minimal number such that 
in the corresponding labelling 
$\{ \gamma (v)\}_{v\in T}$ there is 
$v\in T_n$ with $\gamma (v)\not= 1$. 
Then we claim that $K= Ker\pi_{n}$. 

To see this it is enough to show that 
for every $l\ge n$ every isometry of the form 
$\{ \delta (v)\}_{v\in T}$ such that $T_l$ has 
a unique non-trivial $\delta (w)$ and $w\in layer_l$, 
is contained in $K$. 
Let $g$ and $n$ be as above and $v_0\in T_n$ 
witness that $g$ is of depth $n$. 
We start with the observation that for every 
$l\ge n$ there is $h\in K$ of depth $l$. 
This $h$ can be chosen as an appropriate commutator 
$[g,f]$ where $f\in G$ is of depth $l$.   
For example, if for all $w\in T_l$ which are 
below $v$ we have $\gamma (w)=1$, then we just 
arrange that the labelling of $f$ is non-trivial 
only at two points $w ,w^g \in T_l$ below $v$ 
and the corresponding $(\sigma (w))^{-1}$ is not 
equal to $\sigma (w^g )$.  

To finish the proof of our claim take an isometry 
$g_{\delta}\in K$ defined by $\{ \delta (v)\}_{v\in T}$ 
as in the previous paragraph. 
Choose $w\in layer_l$ with non-trivial $\delta (w)$. 
Conjugating $g_{\delta}$ by elements representing labellings 
with the unique non-trivial element at $w$ we generate 
a subgroup of $K$ having an element $g'$ such that 
the corresponding labelling is non-trivial at $w$ 
but trivial at all points of $T_l \setminus \{ w\}$. 
By sipmlicity of $(P_{l+1} ,A_{l+1} )$, the point 
$w$ in our $g'$ can be labelled by any element of $P_{l+1}$.  
It is also clear that conjugating $g'$ by elements of 
$G$ we can move $w$ to any element of $layer_l$. 
This finishes the proof. 

(2) Here we use implicitly the characterization of 
closed normal subgroups of the group 
$Iso(T)= \wr_{\omega} (Sym (A_i ),A_i )$ from \cite{susz}. 
For every natural $l$ consider the subgroup (denoted by $N_l$ 
of all labellings $\{ \gamma (v)\}_{v\in T}$ such that 
the product $\prod \{ sgn(\gamma (v)):v \in layer_l \}$ is $0$.   
By $sgn(\sigma )$ we denote the parity of the permutation 
$\sigma$ (i.e. it equals to 1 when it is odd). 
It is easy to see that $N_l$ is a normal closed subgroup of  
$Iso(T)$, $l\in \omega$, and $|Iso(T):N_l |=2$. 
Thus the sort of $S(Iso(T))$ consisting of normal open 
subgroups of index 2, is infinite. 
Since $Iso(T)$ is a branch group by Proposition \ref{justinf} 
we see that $S(Iso(T))$ is not $\omega$-categorical. 
Thus it does not satisfy (IP). 
$\Box$ 

\bigskip 

Note that the Iwasawa property is slightly casual 
in the context of branch groups. 

\begin{prop} 
Let $T=T_{\bar{m}}$, where 8 divides some $m_1\cdot ...\cdot m_i$. 
Let $G<Iso(T)$ be a profinite branch group 
with the corresponding system $\{ H_i ,L_i :i\in\omega\}$. 
Then there is a profinite branch group $G^{*}$ which does 
not satisfy (IP), such that for some natural numbers 
$k$ and $k'$ the branch system of $G^{*}$ with indexes 
above $k$ becomes $\{ H_i ,L_i :i\ge k'\}$.  
Moreover $S(G^* )$ is $\omega$-categorical if and 
only if $S(G)$ is $\omega$-categorical. 
\end{prop} 

{\em Proof.} 
To prove this find $k$, such that the number 
$l=|layer_k |$ is of the form $8 l'$.  
Consider the tree $T'=T_{\bar{m}'}$ with $m'_1 =l$ and 
$m'_i = |layer_{k+i-1}|$ for $i\ge 2$. 
The group $G^{*}$ acts on $T'$ as follows. 
Let the kernel $Ker \pi_1 (G^{*} )$ be $H_k$ 
with respect to the action copied from the 
action of $H_k$ on the $k$-th layer of $T$. 
Define $\pi_1 (G^{*})$ to be a regular action 
of $\mathbb{Z}(2)\oplus \mathbb{Z}(4 l')$ on 
the first layer of $T'$.  
Then let $G^{*}$ act on $T'$ as the wreath product of 
these groups. \parskip0pt 

To see that $G^*$ does not satisfy (IP) consider the resulting 
homomorphism $\phi$ from  
$$
G^* \rightarrow \pi_1 (G^*)\rightarrow \mathbb{Z}(2), 
$$ 
where the second one is the projection to $\mathbb{Z}(2)$. 
As an epimorphism $\theta$ of finite groups $H\rightarrow K$ 
take $\mathbb{Z}(4l') \rightarrow \mathbb{Z}(2)$ induced by 
$2l'$-th exponentiation. 
If $\theta:G^* \rightarrow \mathbb{Z}(4l')$ makes this diagram 
commutative, then the generator of $\mathbb{Z}(2)$ in 
the decomposition $\mathbb{Z}(2)\oplus \mathbb{Z}(4l')$  
must have a $2l'$-root. 
This is a contradiction. 

The second statement of the proposition 
follows from the fact that for any profinite group $H$ and 
an open subgroup $P$ the structure $S(H)$ is $\omega$-categorical 
if and only if $S(P)$ is $\omega$-categorical. 
This follows from the fact that the map $K\rightarrow P\cap K$ 
for open subgroups $K<H$ is finite-to-one. 
$\Box$

\bigskip

As we already noted a finitely generated profinite group 
$G$ has finitely many open subgroups of index $n$. 
Thus the structure $S(G)$ for any closed subgroup 
$G<GL_d (\mathbb{Z}_p )$ or $PGL_d (\mathbb{Z}_p )$, 
realizes the simplest case.  
This shows that if any of the groups which we mentioned in 
the introduction, has $\omega$-categorical $S(G)$, then 
it becomes very simple. 
When it is not $\omega$-categorical, then it becomes 
extremely complicated branch group 
(very close to the whole $Iso(T)$). 
This shows that the properties studied in Sections 2 and 3 
correspond to our example much better than $\omega$-categoricity 
of  $S(G)$.


\bigskip

Institute of Mathematics, Wroc{\l}aw University, \parskip0pt

pl.Grunwaldzki 2/4, 50-384 Wroc{\l}aw, Poland \parskip0pt

e-mail: ivanov@math.uni.wroc.pl

\bigskip

\end{document}